\documentclass[reqno]{amsart}
\usepackage{amssymb,latexsym}
\newcommand{\bfi}{\bfseries\itshape}

\newtheorem{thm}{Theorem}[section]

\theoremstyle{definition}

\theoremstyle{remark}
\newtheorem{rem}{Remark}[section]

\makeatletter
\@addtoreset{figure}{section}
\def\thefigure{\thesection.\@arabic\c@figure}
\def\fps@figure{h, t}
\@addtoreset{table}{bsection}
\def\thetable{\thesection.\@arabic\c@table}
\def\fps@table{h, t}
\@addtoreset{equation}{section}

\makeatother

\allowdisplaybreaks

\begin{document}

\title[Averaged Euler equation and a new diffeomorphism group]
{The geometry and analysis of the averaged Euler equations
and a new diffeomorphism group}

\author[J.E. Marsden]{Jerrold E. Marsden}
\address{ Control and Dynamical Systems\\California Institute of
Technology, 107-81\\Pasadena, CA 91125}
\email{marsden@cds.caltech.edu}

\author[T.S. Ratiu]{Tudor S. Ratiu}
\address{D\'{e}partement de Mathematiques\\
Ecole Polytechnique f\'{e}d\'{e}rale de Lausanne\\
CH - 1015 Lausanne, Switzerland}
\email{Tudor.Ratiu@epfl.ch}

\author[S. Shkoller]{Steve Shkoller}
\address{Department of Mathematics\\
University of California \\
Davis, CA 95616}
\email{shkoller@math.ucdavis.edu}


\date{September 1, 1998; current version April 7, 1999}


\begin{abstract}
This paper develops the geometric analysis of geodesic flow
of a new right invariant metric $\langle \cdot, \cdot\rangle_1$ on two
subgroups of the volume preserving
diffeomorphism group of a smooth $n$-dimensional compact subset
$\Omega$ of ${\mathbb R}^n$ with $C^\infty$ boundary $\partial \Omega$.
The geodesic equations are given by the system of PDEs
$$
\begin{array}{c}
\dot{v}(t)+\nabla_{u(t)} v(t)-\epsilon[\nabla u(t)]^t \cdot \triangle u(t)
=  - \operatorname{grad} p(t) \text{ in } \Omega, \\
v=(1-\epsilon \triangle)u, \ \ \operatorname{div} u=0,\\
u(0)=u_0,
\end{array}
$$
which are the averaged Euler (or Euler-$\alpha$) equations
when $\epsilon = \alpha^2$, a length scale, and are the
equations of an inviscid non-newtonian second grade fluid when
$\epsilon = \tilde\alpha_1$, a material parameter.
The boundary conditions associated with the geodesic flow on
the two groups we study are given by either
$$ u =0 \text{ on } \partial \Omega$$
or
$$ u\cdot n =0 \quad \mbox{and} \quad (\nabla_n u)^{{\rm tan}} + S_n( u) =0
\text{ on } \partial \Omega,
$$
where $n$ is the outward pointing unit normal on $\partial \Omega$,
and where $S_n$  is the second
fundamental form of $\partial \Omega$.   We prove that for
initial data $u_0$ in $H^s$, $s >(n/2)+1$,  the above system of
PDE with Dirichlet boundary conditions are well-posed, 
by establishing existence, uniqueness, and
smoothness of the geodesic spray
of the metric $\langle \cdot, \cdot \rangle_1$, together smooth
dependence on initial data.  We are then able to
prove that the limit of zero viscosity
for the corresponding viscous equations is a regular limit.
\end{abstract}

\maketitle


\section{Introduction}
\label{Intro}
\subsection{Background}
The Euler equations of ideal incompressible hydrodynamics on an
$n$-dimensional compact subset $\Omega$ of ${\mathbb R}^n$ with
smooth boundary $\partial \Omega$,
are a  system of partial differential equations describing the motion of a
perfect (ideal, homogeneous, incompressible) fluid and are given by
\begin{equation}\label{euler}
\begin{array}{c}
\partial _t u(t) + \nabla_{u(t)}u(t) = -\text{grad}\ p(t) \text{ in } \Omega\\ u
\text{ is parallel to } \partial \Omega, \\
\operatorname{div} u(t) =0, \ \ u(0)=u_0.
\end{array}
\end{equation}
Here, $p(t)$ is the pressure function which is determined (up
to an additive constant) by the spatial velocity field $u(t)$,
and $\nabla_u u$ denotes the directional (covariant) derivative of $u$ in
the direction $u$; it is often written as $(u \cdot \nabla)u$ using
vector notation in the fluids literature.

The Lagrangian formalism for the hydrodynamics of
incompressible ideal fluids  considers geodesic motion on
$\mathcal{D}_\mu^s:=\mathcal{D}^s_\mu(\Omega)$,
the group of all volume  preserving diffeomorphisms of $\Omega$
of Sobolev class $H^s$.  Geodesics in this context extremize the
energy associated with the $L^2$ norm, which corresponds to the
kinetic energy of the fluid. Arnold \cite{A} and Ebin and Marsden
\cite{EM} showed that $\eta(t)$  is a  smooth geodesic of the weak $L^2$
right invariant metric in $\mathcal{D}^s_\mu$ if and only if the Eulerian
velocity field
$u(t)={\dot{\eta}(t)} \circ \eta(t)^{-1}$
is a solution of the Euler equations.
Moreover, Ebin and Marsden \cite{EM} proved that the geodesic
spray of the $L^2$ right invariant  metric on $\mathcal{D}_\mu^s$
is $C^{\infty}$ for $s>(n/2)+1$.
They derived a number of interesting consequences
from this result, including theorems on the convergence of solutions of the
Navier-Stokes equations to solutions of the Euler equations as the
viscosity limits to zero when $\Omega$ is replaced by a manifold with
no boundary (such as flow in a periodic box).

Marsden, Ebin, and Fischer \cite{MEF} conjectured that although
in a region with boundary, solutions of the Navier-Stokes equations
would {\it not} in general converge to the solutions of the Euler
equations, a certain averaged quantity of the flow may converge.
Recently, Barenblatt and Chorin \cite{BC1,BC2} also speculated that
certain average properties of the flow possess well-defined limits as
the viscosity tends to zero.   This paper proves that an appropriate
choice of right invariant metric on certain subgroups of
$\mathcal{D}_\mu^s$ yields geodesic equations, which may be interpreted
as the ensemble-averaged Euler equations, whose solutions are indeed the
regular limit of the solutions of their  viscous counterparts.

\subsection{Main Results}
We consider two subgroups of $\mathcal{D}^s_\mu$. The first is given
by
$$ \mathcal{D}_{\mu,0}^s = \{\eta \in \mathcal{D}^s_\mu \ | \
\eta = \text{ identity on } {\partial \Omega} \},$$
with $T_e \mathcal{D}_{\mu,0}^s$ consisting of divergence-free $H^s$
vector fields on $\Omega$ that vanish on $\partial\Omega$.

To define the second group,
Let $N$ denote the normal bundle on $\partial \Omega$, and set
$$ {\mathcal N}_\mu^s = \{ \eta \in \mathcal{D}_\mu^s \mid
T\eta|_{\partial \Omega} \cdot n \in H^{s-{\frac{3}{2}}}_\eta(N), \
\;  \mbox{for all} \; \ n \in H^{s-{\frac{1}{2}}}(N)\},$$
where $H^s_\eta$
denotes the space of sections of $N$ covering the diffeomorphism $\eta$.
In the next section, we shall prove that ${\mathcal N}^s_\mu$ is a
$C^\infty$ subgroup of $\mathcal{D}_\mu^s$; see \cite{Shk2} for the
construction of related subgroups of $\mathcal{D}_\mu^s(M)$, for $M$ an
arbitrary compact Riemannian manifold with smooth boundary.

Let ${\mathfrak G}^s_\mu$ denote either $\mathcal{D}_{\mu,0}^s$ or
${\mathcal N}_\mu^s$.
Motivated by the work in \cite{HMR1}, we define a right invariant
$H^1_\alpha$ (pseudo) metric on ${\mathfrak G}_\mu^s$, given at the
identity by
$$
{\frac{1}{2}}\int_\Omega \left[ |u|^2 + \alpha^2 |\nabla u|^2 \right] \mu
+\alpha^2 \int_{\partial \Omega} S_n(u)\cdot u \ \gamma,
$$
where $\alpha >0$ is a constant (representing a length-scale),
$S_n$ is the second fundamental form  of $\partial \Omega$,
and $\gamma$ is the induced ``volume''-form on $\partial \Omega$.

The tangent space of ${\mathcal N}_\mu^s$ at $e$ consists of
divergence-free vector fields of class $H^s$ satisfying the 
free-slip boundary conditions
\begin{equation}\label{fs_bc}
u \cdot n  = 0, \ \ (\nabla_n u)^{{\rm tan}} + S_n(u)=0
\text{ on } {\partial \Omega}.
\end{equation}
Using the Euler-Poincar\'{e} reduction theorem that we recall in
Appendix \ref{app1}, which relates geodesic equations on groups with
their corresponding Euler equations on the associated Lie algebra, we first
show that, formally, geodesics on
${\mathfrak G}^s_\mu$ of the right invariant
$H^1_\alpha$  metric defined above are solutions
of the averaged Euler (or Euler-$\alpha$) equations (see
\cite{HMR1, HMR2}), namely
\begin{equation}\label{ea2}
\begin{array}{c}
\dot{v} + \nabla_u v - \alpha^2 [\nabla u]^T \cdot
\triangle u = -\text{grad }
p \text{ in } \Omega,\\
v=(1-\alpha^2 \triangle) u, \ \ \text{div }u=0,\\
\end{array}
\end{equation}
with either the no-slip boundary conditions
$$ u=0 \text{ on } {\partial \Omega}$$
or the free-slip  boundary conditions
$$ u\cdot n=0, \ \
(\nabla_n u)^{{\rm tan}} + S_n(u)=0
\text{ on } {\partial \Omega}.$$

If the length-scale $\alpha^2$ is replaced
with the material constant $\tilde\alpha_1$, one obtains the equations of
second-grade non-newtonian fluids (see \cite{RE,NT,DF} and references therein).
Notice that the boundary term vanishes if ${\mathfrak G}^s_\mu =
\mathcal{D}_{\mu,0}^s$.

In this paper we shall focus our analysis on the no-slip boundary
condition, as this is the case that has received a great deal of
attention in the literature (see, for example, \cite{CG}, \cite{CO},
and \cite{GGS}).  We shall prove existence and uniqueness of the
geodesic flow of the $H^1_\alpha$ metric on $\mathcal{D}_{\mu,0}^s$
for $s>(n/2)+1$. In fact, we shall prove that the geodesic flow
is $C^\infty$, and has $C^\infty$ dependence on initial data. This
establishes sharp well-posedness on finite time intervals for
classical solutions of the inviscid system (\ref{ea2}).

As a consequence of the smoothness of the geodesic spray on
$\mathcal{D}_{\mu,0}^s$, we are able to prove that
solutions of (\ref{ea2}) with no-slip boundary conditions are a regular
limit of the solutions of the
corresponding viscous equation
\begin{equation}\label{nsa}
\begin{array}{c}
\dot{v} -\nu \triangle u+ \nabla_u v - \alpha^2 [\nabla u]^T \cdot
\triangle u = -\text{grad }
p \text{ in } \Omega,\\
v=(1-\alpha^2 \triangle) u, \ \ \text{div }u=0,\\
u=0 \text{ on } {\partial \Omega},
\end{array}
\end{equation}
answering in the affirmative the conjecture of
Ebin-Fischer-Marsden and Barenblatt-Chorin, as well as establishing the
limit of zero viscosity for second-grade non-newtonian fluids.

The equation (\ref{nsa}) is precisely the equation obtained
from the constitutive theory of simple materials and is the unique
Rivlin-Ericksen momentum equation that satisfies the principles of
material frame indifference and observer objectivity (see \cite{RE, NT,
DF}).  We remark that the mathematical analysis of the viscous equation
(\ref{nsa}) first appeared in the 1984 paper of Cioranescu \& Ouazar
\cite{CO}, where well-posedness on finite time-intervals  for the case of
homogeneous Dirichlet boundary conditions ($u=0$) was established using a
clever eigenfunction expansion for the Galerkin truncation. Using this
technique, Cioranescu and Girault
\cite{CG} were then able to show global existence of (\ref{nsa}) for
small initial data (see also \cite{GGS}). The equations with the stronger
dissipative term $\nu \triangle v $ are studied in
\cite{CF, FHT}.

We mention, finally, that for other problems, such as compressible flow,
the averaged Euler equations and the equations for a non-newtonian fluid
are expected to be different.

\subsection{Outline} The paper is structured as follows.  In Section 2,
we prove that ${\mathcal N}_\mu^s$ is a $C^\infty$ subgroup of
$\mathcal{D}_\mu^s$.
This result uses elliptic operator theory to show that a certain map between
two infinite dimensional vector bundles is a surjection.  In Section 3,
we compute the geodesic spray of the right invariant weak $H^1_\alpha$
(pseudo) metric on $\mathcal{D}_{\mu,0}^s$ and prove that it is a smooth map
in the strong $H^s$ topology for $s>(n/2)+1$.  Finally, in Section 4, we
prove the limit of zero viscosity result.

\section{Subgroups of the diffeomorphism group}
In this section we set up the relevant groups of
diffeomorphisms that we shall need to study the averaged Euler
and second-grade fluid equations in Lagrangian representation.

\subsection{Sobolev Spaces of Mappings.} Let $({M}, g)$
be a compact oriented $C^\infty$ $n$-dimensional
Riemannian manifold  with boundary,
and let $({Q}, g')$ be a
$p$-dimensional compact Riemannian manifold without boundary.
By Sobolev's embedding theorem,
when $s > n/2 + k$, the set of Sobolev mappings $H^s(M,Q)$
 is a subset of $C^k(M,Q)$ with continuous inclusion, and so
for $s>n/2$, an $H^s$-map of $M$ into $Q$ is pointwise
well-defined.  Mappings in the space $H^s(M,Q)$ are those whose
first $s$ distributional derivatives are square integrable in {\it
any} system of charts covering the two manifolds.

For $s>n/2$, the space $H^s(M,Q)$ is a $C^\infty$
differentiable Hilbert manifold.  Let exp$:TQ
\rightarrow Q$ be the exponential mapping associated with $g'$.
Then for each $\phi \in H^s(M,Q)$, the
map $\omega_{\operatorname{exp}}:T_\phi H^s(M,Q) \rightarrow
H^s(M, Q)$ is used to provide a differentiable structure which
is independent of the chosen metric, where
$\omega_{\operatorname{exp}}(v) = \operatorname{exp} \circ v$.

\subsection{Diffeomorphism Groups.} \label{DG}
For a compact Riemannian manifold $M$ with smooth boundary,
the set of $H^s$ mappings from
$M$ to itself is not a smooth manifold; however, if we embed $M$
in its double $\tilde{M}$, then the set $H^s(M,\tilde{M})$ is a
$C^\infty$ Hilbert manifold, and  for $s> n/2+1$, we may form the
set $\mathcal{D}^s(M)$ consisting of $H^s$ maps $\eta$ mapping $M$
to $M$ with $H^s$ inverses. This space {\it is} a smooth manifold.
It is a well-known fact that the diffeomorphism group
$\mathcal{D}^s(M)$ is a
$C^\infty$ topological group for which the left translation
operators are continuous and the right translation operators are
smooth (see \cite{EM} and references therein). One also
knows that $\eta: M
\rightarrow M$ has an extension to an element of (the connected
component of the identity of) $\mathcal{D}^s(\tilde{M})$ if and
only if $\eta$ lies in (the connected component of the identity
of) $\mathcal{D}^s(M)$.

We now restrict our attention to a smooth $n$-dimensional compact subset
$\Omega$ of ${\mathbb R}^n$ with smooth boundary $\partial \Omega$.
Let $\mu=dx^1 \wedge \cdot \cdot \cdot \wedge dx^n$ denote the volume-form
on $\Omega$, and let
$$ \mathcal{D}_\mu^s:= \mathcal{D}_\mu^s(\Omega) := \{\eta \in \mathcal{D}^s
(\Omega) \mid \eta^*(\mu)=\mu \} $$
denote the subgroup of $\mathcal{D}^s(\Omega)$ consisting
of all volume preserving diffeomorphisms of class $H^s$.  For each
$\eta\in \mathcal{D}_\mu^s$, we may use the $L^2$ Hodge
decomposition to define the projection
$P_\eta:T_\eta\mathcal{D}^s \rightarrow T_\eta
\mathcal{D}_\mu^s$ given by
$$ P_\eta (X) = (P_e(X \circ \eta^{-1}))\circ \eta, $$
where $X \in
T_\eta\mathcal{D}_\mu^s$, and  $P_e$ is the $L^2$ orthogonal
projection onto the divergence-free vector fields on $\Omega$. Recall
that this projection is given by
$$ P_e(v) = v - \text{grad }p, $$
where $p$ is the solution of the Neumann problem
$$
\begin{array}{c}
\triangle p = \operatorname{div}\  v  \quad \text{ in } \Omega \\
\frac{\partial p}{\partial n} = v \cdot n \text{ on } \partial \Omega,
\end{array}
$$
and where $n$ is the orientation preserving normal vector field on
$\partial \Omega$.
The function $p$ is the pressure associated with $v$.

\subsection{The subgroup $\mathcal{D}_{\mu,0}^s$}
Ebin \& Marsden \cite{EM} showed that there is a $C^\infty$ differentiable
structure on the those (volume preserving) diffeomorphisms of $\omega$
which keep ${\partial \Omega}$ pointwise fixed.

\begin{thm}\label{thm_em}
The sets
$$ \mathcal{D}_{0}^s = \{ \eta \in \mathcal{D}^s \mid
\eta(x) =x \text{ for all } x \in \partial \Omega\}$$
and
$$ \mathcal{D}_{\mu,0}^s = \{ \eta \in \mathcal{D}_\mu^s \mid
\eta(x) =x \text{ for all } x \in \partial \Omega\}$$
are smooth subgroups of $\mathcal{D}^s$, and $T_e\mathcal{D}^s_0$
consists of $H^s$ vector fields on $\Omega$ vanishing on ${\partial \Omega}$,
while $T_e\mathcal{D}_{\mu,0}^s=\{ u \in T_e\mathcal{D}^s_0|
\operatorname{div} u=0\}$.
\end{thm}
For the proofs, see Section 8 of \cite{EM}.

\subsection{The subgroup ${\mathcal N}_\mu^s$}

For any vector space $E$ and  for all
$\eta \in \mathcal{D}^s$, we set
$H^s_\eta(E):= \{ U \in H^s(\Omega,E) \mid \pi \circ U = \eta\}$, with a similar
definition when  $\eta$ is restricted to $\partial \Omega$.

With $T\Omega = \Omega \times {\mathbb R}^n$,
$$ T\Omega|_{\partial \Omega} = T\partial\Omega \oplus N,$$
where $N$ is the normal bundle.

We define the following vector bundles over $\mathcal{D}_\mu^s$:
\begin{equation}\nonumber
\begin{array}{c}
{\mathcal F} \equiv  \cup_{\eta \in \mathcal{D}_\mu^s}
H^{s-{\frac{3}{2}}}_\eta (T\Omega|\partial \Omega) | \mathcal{D}_\mu^s, \\
{\mathcal E} \equiv  \cup_{\eta \in \mathcal{D}_\mu^s}
H^{s-{\frac{3}{2}}}_\eta (T\partial \Omega) | \mathcal{D}_\mu^s, \\
{\mathcal G} \equiv  \cup_{\eta \in \mathcal{D}_\mu^s}
\left[H^{s-{\frac{3}{2}}}_\eta (T\Omega|\partial \Omega)^*
\otimes H^{s-{\frac{3}{2}}}_\eta (T\partial \Omega)
\right] | \mathcal{D}_\mu^s. \\
\end{array}
\end{equation}

Next, we define the following maps:
\begin{equation}\nonumber
\begin{array}{c}
h: \mathcal{D}_\mu^s \rightarrow {\mathcal F}, \ \
h(\eta) = T\eta|_{\partial \Omega} \cdot n, \ \ n \in H^{s-{\frac{1}{2}}}(N),\\
\Pi:\mathcal{D}_\mu^s \rightarrow {\mathcal G}, \ \
 \Pi(\eta): H^{s-{\frac{3}{2}}}_\eta(T\Omega|\partial \Omega) \rightarrow
 H^{s-{\frac{3}{2}}}_\eta(T\partial \Omega),\\
f: \mathcal{D}_\mu^s \rightarrow {\mathcal E}, \ \ f= \Pi \circ h,
\end{array}
\end{equation}
where $\Pi(\eta)$ is defined pointwise by the ${\mathbb R}^n$-orthogonal
projector $\Pi_{\eta(x)}:T_x\Omega \rightarrow T_x \partial \Omega$ for
$x\in \partial \Omega$.
Lemma B.1 in \cite{Shk2} proves that $f$ is $C^\infty$.

Define the subset ${\mathcal N}_\mu^s$ of $\mathcal{D}_\mu^s$ by
$$ {\mathcal N}_\mu^s = \{ \eta \in \mathcal{D}_\mu^s \mid
T\eta|_{\partial \Omega} \cdot n \in H^{s-{\frac{3}{2}}}_\eta(N), \
\;  \mbox{for all} \; \ n \in H^{s-{\frac{1}{2}}}(N)\}.$$

\begin{thm}\label{thm1}
The set ${\mathcal N}_\mu^s$ is a subgroup of $\mathcal{D}_\mu^s$
for $s>{\frac{n}{2}}+1$, such that
\begin{eqnarray*}
&&T_e {\mathcal N}_\mu^s = \{ u \in T_e \mathcal{D}_\mu^s \mid
(\nabla_n u)^{\rm tan} + S_n(u) =0 \text{ on } \partial \Omega
\ \ \forall n \in H^{s-{\frac{1}{2}}}(N) \},
\end{eqnarray*}
where $S_n:T\partial \Omega \rightarrow T\partial \Omega$
is the second fundamental form of $\partial \Omega$ given by
$$\langle S_n \langle u \rangle, v \rangle =-\langle \nabla_un,v\rangle, \ \
u,v \in H^{s-{\frac{3}{2}}}(T{\partial \Omega}).$$
\end{thm}

\begin{proof}
It is clear that ${\mathcal N}_\mu^s$ is closed under right composition;
hence, we must show that ${\mathcal N}_\mu^s$
is a submanifold of $\mathcal{D}_\mu^s$.  To do so, we shall use the
transversal mapping theorem (see, for example, \cite{AMR}) which
states that if $f:\mathcal{D}_\mu^s \rightarrow {\mathcal E}$ is
transversal to the zero section of ${\mathcal E}$, then ${\mathcal N}_\mu^s
= f^{-1}(0)$ is a submanifold of $\mathcal{D}_\mu^s$.

Since our manifolds are Hilbert, in order to establish the transversality of
$f$ with $0 \in C^\infty({\mathcal E})$, it suffices to prove that
$f$ is a surjection.  The Frechet derivative on ${\mathbb R}^n$
induces, by a pointwise lift, natural (weak) covariant derivatives
$\overline \nabla$ on ${\mathcal F}$ and ${\mathcal G}$
(see Lemma B.1 in \cite{Shk2} and Section 9 of \cite{EM}).

We compute that for all $u$ in $T_\eta\mathcal{D}_\mu^s
= H^s_\eta(T\Omega)$,
\begin{equation}\label{h_der}
\overline\nabla_u h(\eta) = \nabla_n u,
\end{equation}
where $\nabla$ denotes
the covariant derivative in the pull-back bundle $\eta^*(T\Omega)$

Next, we compute the covariant derivative of $\Pi$.  For all
$u \in T_\eta \mathcal{D}_\mu^s$, and $v,z \in {\mathcal F}_\eta$, along
the boundary $\partial \Omega$,
\begin{eqnarray}
&& [\nabla _u \Pi_{\eta(x)}] (v(x)) \cdot z(x)
= - (\nabla_u v(x))^{{\rm tan}} \cdot z(x) \nonumber\\
&& \qquad \qquad  -  (\nabla_u z(x))^{{\rm tan}} \cdot
v(x) - u\cdot \nabla \left[ v^{{\rm tan}} (x) \cdot z^{{\rm tan}}(x)
\right], \label{1}
\end{eqnarray}
where
$(\cdot)^{{\rm tan}}$ denotes the tangential component.  Hence, $\nabla _u
\Pi_\eta$
is symmetric with respect to the inner-product on ${\mathbb R}^n$.  Now,
by definition, for $x \in {\partial \Omega}$,
$$\Pi_{\eta(x)}(v(x)) \cdot \nu(x) = 0 \
\;  \mbox{for all} \; \ v \in {\mathcal F}_\eta,  \nu \in
H^{s-{\frac{3}{2}}}_\eta(N),$$
so setting $v = \nu$ in (\ref{1}) shows that
\begin{equation}\label{P_der}
[\overline\nabla _u \Pi_\eta] (\nu) = -(\nabla_u \nu)^{{\rm tan}} =
S_\nu (u).
\end{equation}

It follows that  for all $\eta \in f^{-1}(0)$,
\begin{align*}
\overline\nabla_u f(\eta) &=
\overline\nabla_u \Pi_\eta \cdot h(\eta) + \Pi_\eta \cdot
\overline\nabla_u h(\eta) \\
&=  S_\nu (u)  + (\nabla_n u)^{\rm tan} \in {\mathcal E}_\eta,
\end{align*}
where $\nu = T\eta \cdot n \in H^{s-\frac{3}{2}}_\eta(N)$.

It remains to show that for every $w \in {\mathcal E}_\eta$, there
exists $u\in T_\eta\mathcal{D}_\mu^s$ such that $\overline\nabla_u f(\eta)
 =w$.
By right translation to the identity, it suffices to find
$u \in T_e\mathcal{D}_\mu^s$ such that $\overline\nabla_u f(e)= w$ for every
$w \in H^{s -{\frac{3}{2}}}(T\partial \Omega)$.

To do so we obtain a solution to the following elliptic boundary
value problem:
For $F \in H^{s-2}(T\Omega)$, $w \in H^{s -{\frac{3}{2}}}(T\partial
\Omega)$ and
$n \in H^{s-{\frac{1}{2}}}(N)$, find
$(u,p) \in T_e\mathcal{D}_\mu^s \times H^{s-1}(\Omega)/{\mathbb R}$ such that
\begin{equation}\label{stokes1}
\begin{array}{c}
-\triangle u + \text{grad }p= F \text{ in } \Omega\\
(\nabla_n u)^{\rm tan}+ S^n(u) =w \text{ on } \partial \Omega,
\end{array}
\end{equation}
where in Cartesian components $(\triangle u)^i = \partial_k\partial_k u_i $,
i.e. $\triangle$ is the component-wise Laplace operator.
Note that by definition of $T_e\mathcal{D}_\mu^s$,
$u\cdot n =0$ on ${\partial \Omega}$ and $\operatorname{div}u =0$.

A weak solution to (\ref{stokes1}) in the class of $H^1$ divergence-free
vector fields that are parallel to $\partial \Omega$ is supplied by Step 3
of the proof of
Theorem 2.1 in \cite{Shk2}.  Noting that in coordinates
$(\nabla_nu)^{\rm tan} + S_n(u) = [(u_{i,j} + u_{j,i})\cdot n_j]^{\rm tan}$,
Theorem 2.8 of \cite{MPS} provides a strong solution in the class of
$H^2$ divergence-free vector fields that are parallel to the boundary
and satisfy the boundary condition $(\nabla_nu)^{\rm tan} + S_n(u) =0$
on $\partial \Omega$.
Theorem 2.9 of \cite{MPS} then provides the
elliptic regularity required to obtain $u \in T_e\mathcal{D}_\mu^s$
whenever
$(F,w) \in H^{s-2}(T\Omega) \times H^{s-{\frac{3}{2}}}(T{\partial \Omega})$,
and this completes the proof of the theorem.
\end{proof}

A similar argument also yields

\begin{thm}\label{thm1a}
The set
$$ {\mathcal N}^s = \{ \eta \in \mathcal{D}^s \mid
T\eta|{\partial \Omega} \cdot n \in H^{s-{\frac{3}{2}}}_\eta(N), \
\;  \mbox{for all} \; \ n \in H^{s-{\frac{1}{2}}}(N)\}$$
is a subgroup of $\mathcal{D}^s$ for $s>{\frac{n}{2}}+1$, such that
\begin{eqnarray*}
&&T_e {\mathcal N}^s = \{ u \in T_e \mathcal{D}^s \mid
(\nabla_n u|_{\partial \Omega})^{{\rm tan}}
+ S_n(u)=0 \text{ on } \partial\Omega
\ \ \forall n \in H^{s-{\frac{1}{2}}}(N) \}.
\end{eqnarray*}
\end{thm}
See \cite{Shk2} for the construction of $C^\infty$ differentiable structure on
a number of new diffeomorphism groups of arbitrary compact Riemannian manifolds
with boundary that
describe particular hydrodynamic motions.

Now let ${\mathfrak G}^s =\mathcal{D}_{0}^s$  or ${\mathcal N}^s$, and
let ${\mathfrak G}^s_\mu =\mathcal{D}_{\mu,0}^s$  or ${\mathcal N}_\mu^s$.

We do not call $T_e{\mathfrak G}_\mu^s$ literally the Lie algebra
of ${\mathfrak G}_\mu^s$ as the bracket losses regularity and thus
does not belong to the Hilbert class $H^s$; nevertheless, the bracket
$[u,v]$ of two elements $u, v \in T_e{\mathfrak G}_\mu^s$
satisfies the boundary conditions (\ref{fs_bc}).  To see this, let
$\psi_t$ be the flow of $u$ and $\phi_t$ the flow of $v$.  Then
the flow $\sigma_t$ of $[u,v]$ may be expressed as
$$\sigma_t = \lim_{n\rightarrow \infty} (
\phi_{-\sqrt{t/n}} \circ \psi_{-\sqrt{t/n}} \circ
\phi_{\sqrt{t/n}} \circ \psi_{\sqrt{t/n}})^n,$$
for $t \ge 0$.  Hence, it is clear that
$T{\sigma_t}|_{\partial \Omega}$ maps sections of $N$ into sections of $N$
so that $[u,v]$ must satisfy (\ref{fs_bc}).

\subsection{The projector ${\mathcal P}$}
For $r\ge 1$, let  ${\mathcal V}^r$ denote the $H^r$ vector fields on
$\Omega$ which satisfy the boundary conditions prescribed  to elements of
$T_e{\mathfrak G}^s$ and let ${\mathcal V}^r_\mu = \{ u\in {\mathcal V}^r \
| \ \text{div }u=0\}$.

Define the Stokes projector  by
\begin{equation}\label{P}
\begin{array}{c}
{\mathcal P}_e:{\mathcal V}^r  \rightarrow {\mathcal V}^r_\mu,\\
{\mathcal P}(w) = w -(1-\triangle)^{-1}\text{grad }p,
\end{array}
\end{equation}
where $p$ depends on $v$ and the pair $(v,p) \in {\mathcal V}^r_\mu
\times H^{r-1}(\Omega)/{\mathbb R}$ solves the Stokes problem
\begin{equation}\nonumber
\begin{array}{c}
(1-\triangle) v + \text{grad }p = (1-\triangle) w,\\
\text{div } v =0,
\end{array}
\end{equation}
together with either the no-slip or free-slip boundary conditions, as
appropriate. The Stokes projector ${\mathcal P}_e$ induces the
decomposition
$${\mathcal V}^r={\mathcal V}^r_\mu \oplus (1-\triangle)^{-1}\text{grad}
H^{r-1}(\Omega),$$
and it is readily checked that the two summands are orthogonal with
respect to $\langle \cdot, \cdot \rangle_1(e)$.

For $s>(n/2)+1$, define ${\mathcal P}:T{\mathfrak G}^s \rightarrow
T{\mathfrak G}^s_\mu$ to be the bundle map covering the identity, given
on each fiber by
\begin{equation}\nonumber
\begin{array}{c}
{\mathcal P}_\eta:T_\eta{\mathfrak G}^s \rightarrow T_\eta{\mathfrak
G}^s_\mu,\\
{\mathcal P}_\eta(X_\eta)= \left[ {\mathcal P}_e (X_\eta \circ \eta^{-1})
\right] \circ \eta.
\end{array}
\end{equation}
Theorem 3.1 in \cite{Shk2} proves that ${\mathcal P}$ is a well-defined
$C^\infty$ bundle map; this fact will be crucial in proving that the 
geodesic spray of
the invariant metric $\langle\cdot,\cdot\rangle_1$ is smooth.

\section{Geodesic motion}
\label{sec_spray}

Again, let ${\mathfrak G}^s =\mathcal{D}_{0}^s$  or ${\mathcal N}^s$, and
let ${\mathfrak G}^s_\mu =\mathcal{D}_{\mu,0}^s$  or ${\mathcal N}_\mu^s$.

\subsection{$H^1_\alpha$ metric on ${\mathfrak G}_\mu^s$}
In this section, we shall analyse the geodesic motion of the weak $H^1_\alpha$
right invariant (pseudo) metric $\langle \cdot, \cdot\rangle_1$ on  the group
${\mathfrak G}_\mu^s$.  This metric is defined as follows:
For $X,Y \in T_e {\mathfrak G}_\mu^s$, we set
\begin{align}
\langle X, Y \rangle_{1}(e) &=  \int_\Omega \left( X(x)\cdot Y(x) +
\alpha^2 \nabla X(x)\cdot\nabla Y(x)\right) \mu(x)
\nonumber \\
&  \qquad \qquad +
 \alpha^2 \int_{\partial \Omega} S_n(X(x))\cdot Y(x)\ \gamma(x)
\label{s1}
\end{align}
and extend $\langle \cdot, \cdot \rangle_{1}$ to ${\mathfrak G}_\mu^s$
by right invariance.  Here $n$ is the outward unit normal on
${\partial \Omega}$ and $\gamma$ is the induced volume measure on
${\partial \Omega}$.  Again, if ${\mathfrak G}^s_\mu=\mathcal{D}_{\mu,0}^s$,
then the boundary term vanishes.

\subsection{Euler-Poincar\'{e} equations on $T_e{\mathfrak G}_\mu^s$}
Appendix \ref{app1} is devoted to a review of Lagrangian reduction on
topological groups with a one-sided invariant metric which leads to a system
of reduced equations that are called the Euler-Poincar\'{e} equations.
We refer the reader to this appendix for the general theory; for
purposes of the current development, we shall restrict attention to
geodesics of a right invariant metric on ${\mathfrak G}_\mu^s$.
The fundamental idea is to use the $C^\infty$ right translation maps
on ${\mathfrak G}_\mu^s$ to translate geodesic
motion over the entire topological group ${\mathfrak G}_\mu^s$ onto
motion in the single fiber $T_e{\mathfrak G}_\mu^s$. We shall state
the reduction theorem in this context.

\begin{thm}[Euler-Poincar\'{e} for ${\mathfrak G}_\mu^s$]\label{thm_ep}
Consider ${\mathfrak G}_\mu^s$ with the right invariant metric
$\langle \cdot, \cdot\rangle_{1}$.
A curve $\eta(t)$ in
${\mathfrak G}_\mu^s$ is a geodesic of this metric if and only if
$u(t) = T_{\eta(t)}R_{\eta(t)^{-1}} \dot \eta(t) = {\dot{\eta}}(t) \circ
\eta(t) ^{-1}$
satisfies
\begin{equation}\label{EP}
\frac{d}{dt} u(t) = - \operatorname{ad}^*_{u(t)}u(t)
\end{equation}
where $\operatorname{ad}^*_u$ is the formal adjoint of
$\operatorname{ad}_u$ with respect to the inner-product $\langle \cdot,
\cdot \rangle_1(e)$ given by
$$ \langle \operatorname{ad}^*_uv,w \rangle_{1}(e) = \langle v, [u,w]
\rangle_1(e)$$
for all $ u,v,w \in T_e {\mathfrak G}_\mu^s$, where
\begin{equation}\nonumber
\operatorname{ad}^*_uu = (1- \alpha^2\triangle)^{-1}\left[
\nabla_{u(t)} (1-\alpha^2\triangle)u(t)
- \alpha^2[\nabla u(t)]^t \cdot \triangle u(t) +
 \operatorname{grad} p(t)\right],
\end{equation}
and $u=0$ on $\partial \Omega$ if ${\mathfrak G}^s_\mu =
\mathcal{D}_{\mu,0}^s$ and $u\cdot n=0$, $(\nabla_nu)^{\rm tan}
+S_n(u)=0$ on $\partial \Omega$ if ${\mathfrak G}^s_\mu={\mathcal N}_\mu^s$.
\end{thm}
\begin{proof}
Restricting $\langle \cdot, \cdot\rangle_1$ to the algebra
$T_e \mathfrak{G}_\mu^s$, we compute the first variation of the
action function
$$ {\frac{1}{2}}\int _a ^b\int_\Omega \left[|u|^2
+\alpha^2  |\nabla u|^2 \right] \mu dt
+ \alpha^2 \int_{\partial \Omega} S_n(u)\cdot u \ \gamma dt$$
for constrained variations of the form $\delta u = \partial_t w - [w,u]$.
Integrating by parts, we obtain that
$$\int _a ^b \int_\Omega \left[(1-\alpha^2\triangle)u\right] \cdot
\left[ \partial_t w + [w,u]\right]  \mu dt +
\int _a ^b\int_{\partial \Omega}
 \alpha^2 \left[ \nabla_n u \cdot \delta u +
S_n(u)\cdot \delta u \right] \gamma dt.$$

The boundary term vanishes for $u$ and $\delta u$ in $T_e{\mathfrak G}^s_\mu$,
so another integration by parts yields
\begin{align*}
\int _a ^b \int_\Omega &\left[
\partial_t(1-\alpha^2\triangle)u + \nabla_u(1-\alpha^2\triangle)u -\alpha^2
\nabla u^t \cdot \triangle u \right] \cdot w \ \mu dt \\
&=
\int _a ^b \langle \partial_tu +
(1-\alpha^2\triangle)^{-1}\left[\nabla_u(1-\alpha^2\triangle)u -\alpha^2
\nabla u^t \cdot \triangle u \right], w \rangle_1(e) dt.
\end{align*}
Since $w \in T_e{\mathfrak G}^s_\mu$ is arbitrary, $u$ is a fixed-point
of the action if and only if
$$\partial_tu+ {\mathcal P}_e \circ
(1-\alpha^2\triangle)^{-1}\left[\nabla_u(1-\alpha^2\triangle)u -\alpha^2
\nabla u^t \cdot \triangle u \right] =0.$$
Using the definition of the Stokes projector ${\mathcal P}_e$ concludes
the proof.
\end{proof}

The Euler--Poincar\'e equation
\begin{equation}\label{ss}
\begin{array}{c}
\partial_t (1-\alpha^2\triangle)u + \nabla_u (1-\alpha^2\triangle)u
- \alpha^2[\nabla u)]^t \cdot \triangle u= - \text{grad } p,\\
\text{div }u=0, \ \ u(0)=u_0,
\end{array}
\end{equation}
together with either the no-slip boundary condition $u=0$ or the free-slip
boundary condition $u\cdot n=0$ and $(\nabla_n u)^{\rm tan} + S_n(u)=0$,
is called the {\bfi averaged Euler equation\/} or
the {\bfi Euler-$\alpha$ equation\/}.  As we already mentioned, this equation
is also the equation for {\bfi inviscid second-grade non-Newtonian fluids} when
$\alpha^2$ is replaced by $\tilde \alpha_1$, a material parameter measuring
the elastic response of the fluid.

Being Euler--Poincar\'e equations, of course these equations share
all the properties given by the usual Euler equations, such as a
Kelvin-Noether theorem, a Lie--Poisson Hamiltonian structure and
so on (see \cite{HMR2} for some of the basic facts and literature).

\subsection{The geodesic equations on $T^*_e{\mathfrak G}_\mu^s$}

On the dual of $T_e{\mathfrak G}_\mu^s$, a simple
computation of the coadjoint action verifies that the averaged Euler or
 Euler-$\alpha$ equations may be expressed as
$$ \partial_t v + {\mathcal L}_u v = -dp,$$
where the one-form $v$ is associated to $u^\flat$ by
$v=(1-\alpha^2 \triangle) u^\flat$.  Using the exterior derivative $d$,
we may identify the dual of $T_e{\mathfrak G}_\mu^s$ with two-forms $\omega$
(as in \cite{MW})
and write the Euler-$\alpha$ equations in vorticity form as
$$ \partial_t \omega + {\mathcal L}^\alpha_u \omega =0,$$
where $\omega = d u^\flat$, and\footnote{
More generally, geodesics of the $H^m$ right invariant metric are
defined as above, but with the conjugated Lie derivative operator
${\mathcal L}^{m,\alpha}_u = (1-\alpha^2 \triangle)^{-m} {\mathcal L}_u
(1-\alpha^2 \triangle)^m.$}
$$ {\mathcal L}^\alpha_u = (1-\alpha^2 \triangle)^{-1} {\mathcal L}_u
(1-\alpha^2 \triangle).$$

For example, on ${\mathbb T}^2$, we identify $T^*_e{\mathfrak G}_\mu^s$
with smooth functions, and write the Euler-$\alpha$ equations as
$$ \partial_t q + \nabla_u q=0, \ \ q=(1-\alpha^2\triangle)\omega,$$
where $\omega = d u^\flat$.  Letting $\omega = -\triangle \psi$,
these equations take the familiar Lie-Poisson form
$$ \partial_t q = \{ \psi, q\}.$$

\subsection{Smoothness of the geodesic spray on $\mathcal{D}_{\mu^s}$}
In this section, we establish the well-posedness of the averaged Euler
equations with no-slip boundary conditions by proving that the
geodesic spray of $\langle\cdot,\cdot\rangle_1$ is smooth on
$\mathcal{D}_{\mu,0}^s$.  (See Theorem 3.3 of \cite{Shk1} for the
smoothness of the geodesic spray on $\mathcal{D}_\mu^s(M)$ when
$M$ is an arbitrary compact boundaryless Riemannian manifold.)

\begin{thm}
\label{thm_spray1}
For $s>{\frac{n}{2}}+1$ and $u_0 \in
\mathcal{D}_{\mu,0}^s$, there exists an open interval $I=(-T,T)$ depending on
on $u_0$, and a unique geodesic $\dot\eta$ of $\langle\cdot,\cdot\rangle_1$
such
that $u=\dot\eta \circ \eta$ satisfies (\ref{ss}) with $\eta(0)=e$
and $\dot\eta(0)=u_0$ such that
$\dot\eta \in C^\infty(I,T\mathcal{D}_{\mu,0}^s)$ has $C^\infty$
dependence on $u_0$
\end{thm}
\begin{proof}
We compute the first variation of the
action function
$$
{\mathcal E}(\eta) = {\frac{1}{2}}\int _a ^b \langle
{\dot{\eta}}(t),  {\dot{\eta}}(t) \rangle_{1} dt,
$$
which we decompose as
$$
{\mathcal E}^0(\eta) = \frac{1}{2} \int _a ^b \int_\Omega
|\dot{\eta}(x)|^2\mu(x) dt
$$
and
$$
{\mathcal E}^1(\eta)= \frac{\alpha^2}{2 }\int _a ^b
\int_\Omega
|\nabla (\dot{\eta}\circ\eta^{-1})(y)|^2 \mu(y) dt.
$$
By definition of $\mathcal{D}_{\mu,0}^s$, the boundary terms appearing from
integration by parts  vanish; hence, we restrict our computations to
the interior.
We have
\begin{eqnarray*}
{\mathcal E}^1(\eta)
&=& {\frac{\alpha^2}{2}}
\int _a ^b \int_\Omega |
\nabla({\dot{\eta}}\circ\eta^{-1})(y)|^2 dy dt \\
&=& {\frac{\alpha^2}{2}} \int _a ^b
\int_\Omega \left( \nabla {\dot{\eta}}(x)
\cdot [T\eta(x)]^{-1}\right) \cdot\left(\nabla {\dot{\eta}}(x)
\cdot [T\eta(x)]^{-1}\right) dxdt.
\end{eqnarray*}
Let $\epsilon \mapsto \eta^\epsilon$ be a smooth curve in
$\mathcal{D}_{\mu,0}^s$ such that $\eta^0=\eta$ and
$(d /d\epsilon)|_{ \epsilon = 0} \eta^\epsilon = \delta
\eta$.  Then
\begin{align*}
\mathbf{d} {\mathcal E}^1(\eta)\cdot \delta \eta
&=  \alpha^2\int _a ^b\int_\Omega
\left( \left. \frac{D}{d\epsilon} \right|_0
\nabla {\dot{\eta}}^\epsilon \cdot [T\eta^\epsilon]^{-1} \right)
\cdot \left( \nabla {\dot{\eta}} \cdot [T\eta]^{-1} \right) dx dt\\
&=  \alpha^2\int _a ^b\int_\Omega
\left[
\left(
\left. \frac{D}{d\epsilon} \right|_0
      \nabla{\dot{\eta}}^\epsilon \cdot [T\eta]^{-1} \right)
\cdot \bigl(
\nabla{\dot\eta} \cdot [T\eta]^{-1} \bigr) \right.  \\
& \quad \left. - \left( \nabla\delta\eta \cdot  [T\eta]^{-1}\right) \cdot
\left( \nabla{\dot{\eta}} \cdot [T\eta]^{-1})^t
(\nabla{\dot{\eta}} \cdot [T\eta]^{-1}\right)
\right] dx dt \\
&=  \alpha^2\int _a ^b\int_\Omega \left[
\left( \nabla [ (D/dt) \delta \eta] \cdot [T\eta]^{-1} \right) \cdot
\left( \nabla\dot\eta \cdot [T\eta]^{-1} \right) \right.\\
&  \quad -
\left. \left( \nabla \delta \eta\right) \cdot  \left( (\nabla {\dot{\eta}}
\cdot
[T\eta]^{-1})^t
\cdot (\nabla {\dot{\eta}}\cdot [T\eta]^{-1}) \cdot {[T\eta]^{-1}}^t \right)
\right] dx dt.
\end{align*}
Integration by parts yields
\begin{eqnarray*}
&& \int _a ^b\int_\Omega
\left(\nabla (\frac{D}{dt} \delta \eta) \cdot [T\eta]^{-1}\right) \cdot
\left( \nabla {\dot{\eta}} \cdot[T\eta]^{-1} \right) dx dt \\
&&\qquad \qquad =
-\int _a ^b\int_\Omega
\left( \nabla \delta \eta\right) \cdot
\left( \frac{D}{dt} \{ \nabla {\dot{\eta}} \cdot
[T\eta]^{-1}\cdot {[T\eta]^{-1}}^t\}\right) dx dt
\end{eqnarray*}
We use the product rule to get that
\begin{align*}
\frac{D}{dt} \{ \nabla {\dot{\eta}}
\cdot [T\eta]^{-1} \cdot {[T\eta]^{-1}}^t\}
&=
\nabla {\ddot{\eta}}
\cdot [T\eta]^{-1} \cdot {[T\eta]^{-1}}^t \\
&  - ( \nabla{\dot{\eta}}
\cdot [T\eta]^{-1}) \cdot (\nabla {\dot{\eta}} \cdot
[T\eta]^{-1}) \cdot {[T\eta]^{-1}}^t\\
&  -(\nabla{\dot{\eta}}\cdot [T\eta]^{-1})
\cdot (\nabla {\dot{\eta}} \cdot [T\eta]^{-1})^t
\cdot {[T\eta]^{-1}}^t.
\end{align*}
Integrating by parts, noting that the boundary terms vanish by
virtue of the subgroup
$\mathcal{D}_{\mu,0}^s$,  we have that
\begin{align*}
\mathbf{d} {\mathcal E}^1(\eta) \cdot \delta \eta
= & \; \alpha^2 \int _a ^b\int_\Omega \bigl[
\text{div}\bigl( \bigl\{ \nabla \ddot\eta \cdot [T\eta]^{-1}
-(\nabla \dot\eta \cdot [T\eta]^{-1})^t\cdot (\nabla \dot\eta \cdot
[T\eta]^{-1}) \\
& \qquad \qquad
+(\nabla \dot\eta \cdot [T\eta]^{-1})\cdot (\nabla \dot\eta \cdot
[T\eta]^{-1}) \bigr\}
\cdot {[T\eta]^{-1}}^t
\bigr)
\end{align*}

Computing the first variation of ${\mathcal E}^0$, we obtain
\begin{align*}
\mathbf{d} {\mathcal E}^0(\eta) \cdot \delta \eta &=
\int _a ^b\int_\Omega
\left.\frac{D}{d\epsilon}\right|_0 {\dot{\eta}}^\epsilon \cdot
{\dot{\eta}} \ dx dt =
\int _a ^b\int_\Omega
\frac{D}{dt} \delta\eta \cdot
{\dot{\eta}}\ dx dt \\
&=
\int _a ^b\int_\Omega
 -{\ddot{\eta}}\cdot
 \delta \eta \ dx dt .
\end{align*}
Setting $\mathbf{d} {\mathcal E} \cdot \delta \eta = 0$, and
using the projector
${\mathcal P}$ given by (\ref{P}) gives
\begin{align*}
{\mathcal P}_\eta \circ {\ddot{\eta}}&=
{\mathcal P}_\eta \circ (1- \alpha^2\widehat{\triangle}_\eta)^{-1}
\left[
\text{div}\left\{ (\nabla{\dot{\eta}}
[T\eta]^{-1})^t (\nabla {\dot{\eta}}
[T\eta]^{-1}) \right.\right.\\
&\qquad \left. - \nabla{\dot{\eta}} [T\eta]^{-1}
\nabla{\dot{\eta}} [T\eta]^{-1}
 - (\nabla{\dot{\eta}}[T\eta]^{-1})
( \nabla{\dot{\eta}}[T\eta]^{-1})^t
\right\} {[T\eta]^{-1}}^t\bigr],
 \label{spray}
\end{align*}
where
\begin{equation}
\label{hatlap.equation}
\widehat{\triangle}_\eta=\text{div}[\nabla(\cdot )
(T\eta)^{-1} {(T\eta)^{-1}}^t].
\end{equation}

Let us prove that the above expression is well-defined; namely, we shall
show that it makes for the Stokes projector to act on both $\ddot\eta$
and $F_\eta$.  To see this, notice that $\widehat{\triangle}_\eta=
TR_\eta \circ \triangle \circ TR_{\eta^{-1}}$, and  that
${\mathcal P}_\eta(\ddot \eta) = [{\mathcal P}_e(\partial_t + \nabla_uu)]
\circ \eta$, where $u=\dot\eta\circ\eta^{-1}$.
The Stokes operator acts on $(\partial_t + \nabla_uu)$ by  $(1-\alpha^2
\triangle)$ whose domain is $H^2(T\Omega) \cap H^1_0(T\Omega)$, and this
operation is well-defined as both $\partial_t u$ and $\nabla_uu$ are
in the domain of $(1-\alpha^2\triangle)$, since $u=0$ on $\partial \Omega$.

We may reexpress the above equation as
\begin{align*}
[\partial_t u + {\mathcal P}_e(\nabla_uu)]\circ \eta =
[{\mathcal P}_e \circ (1-\alpha^2\triangle)^{-1}\text{div}[
\nabla u^t\cdot \nabla u - \nabla u \cdot \nabla u - \nabla u\cdot
\nabla u^t]\circ \eta;
\end{align*}
thus the right-hand-side is also well-defined as the image of $(1-\alpha^2
\triangle)^{-1}$ is the domain of $(1-\alpha^2 \triangle)$.  Denoting
the right-hand-side of the above equation by ${\mathcal S}_\eta(\dot\eta)$,
we have that
$$\ddot \eta = {\mathcal B}(\eta,\dot\eta):=
(1-{\mathcal P}_\eta)\circ \ddot \eta +
{\mathcal S}_\eta(\dot\eta).$$
We rewrite this equation as the system
\begin{align*}
\dot\eta&= V_\eta,\\
\ddot \eta& =\frac{dV_\eta}{dt} = {\mathcal B}(\eta,\dot\eta),\\
\eta&(0)=e,\ \ V_\eta(0)=u_0.
\end{align*}
We shall prove that ${\mathcal B}:T\mathcal{D}_{\mu,0}^s \rightarrow
T^2\mathcal{D}_{\mu,0}^s$ and that ${\mathcal B}$ is a $C^\infty$ bundle
map.  Then the
standard theorem for existence and uniqueness of ordinary
differential equations on a Hilbert manifold provides the existence
of a unique $C^\infty$ curve $\dot\eta(t)$ solving the above system
on $[0,T)$, that depends smoothly on the initial data $u_0$;
the time-reversal symmetry allows us to extend the interval to $(-T,T)$.

That ${\mathcal B}$ is $C^\infty$ follows from the fact that
$\nabla u \cdot \nabla u$ is of class $H^{s-1}$ whenever
$u$ is in $H^s$ (because $H^{s-1}$ forms
a multiplicative algebra when $s>(n/2)+1$), so that
$(1-\alpha^2\triangle)^{-1}\text{div}[
\nabla u^t\cdot \nabla u - \nabla u \cdot \nabla u - \nabla u\cdot
\nabla u^t]$ is in $H^s(T\Omega) \cap H^1_0(T\Omega)$.
That $S_\eta(\dot\eta)$ is of class $H^s$ follows from Theorem B.1 in
\cite{Shk2} together with the smoothness of the Stokes projector.

The fact that $(1-{\mathcal P}_\eta)\circ \ddot \eta$ is of class $H^s$
whenever $\dot\eta \in T_\eta\mathcal{D}_{\mu,0}^s$ follows from
similar arguments (see the proof of Theorem 4.1 in \cite{Shk2} for
details).
\end{proof}

Using the fact that the inversion map $\eta \mapsto \eta^{-1}$ is
only $C^0$ as a map of $\mathcal{D}_{\mu,0}^s$ into
$\mathcal{D}_{\mu,0}^s$, and is $C^1$ as a map of
$\mathcal{D}_{\mu,0}^s$ into
$\mathcal{D}_{\mu,0}^{s-1}$, we immediately obtain that
$$ u \in C^0( I, {\mathcal V}^s_\mu) \cap C^1(I,{\mathcal V}^{s-1}_\mu)$$
where for $r\ge 1$, ${\mathcal V}^r_\mu =\{ u\in H^s(T\Omega)\cap H^1_0
(T\Omega) \ | \ \text{div }u=0\}$, and has $C^0$ dependence on the initial
data $u_0$.

\section{The regular limit of zero viscosity}

The viscous averaged Euler equations also termed the averaged Navier-Stokes
equations  or the viscous equations of second-grade non-newtonian
fluids are given by
\begin{equation}\label{nse}
\begin{array}{c}
\partial_t (1-\alpha^2\triangle)u - \nu \triangle u
+ \left[ \nabla_u (1-\alpha^2\triangle)u
- \alpha^2 \nabla u^t \cdot\triangle u\right]=
 - \operatorname{grad} \ p,\\
\text{div }u=0,\\
u=0 \text{ on } \partial\Omega, \ \ u(0)=u_0.
\end{array}
\end{equation}
In \cite{CG}, well-posedness of (\ref{nse}) was established, but
in 3D, the estimates relied crucially on the presense of viscosity,
so that a limit of zero viscosity theorem did not follow.
Having proven the smoothness of the geodesic
spray of the Euler-$\alpha$ equations, we follow \cite{EM}
and use the product formula approach to prove the existence of viscosity
independent solutions to (\ref{nse}) on finite time intervals as well
as the existence of the limit of zero viscosity.  In the case that
$\alpha=0$, this limiting procedure is believed to be valid only
for compact manifolds without boundary (e.g., for flows with
periodic boundary conditions), as the Navier-Stokes equations and
the Euler equations do not share the same boundary conditions on
manifolds with boundary.

\begin{thm}\label{product}
Let ${\mathcal B}: T\mathcal{D}_{\mu,0}^s\rightarrow T^2\mathcal{D}_{\mu,0}^s$
be the $C^\infty$ geodesic spray of the metric $\langle\cdot,\cdot\rangle_1$.
For each $s>(n/2)+1$, let ${\mathcal T}:T_e \mathcal{D}_{\mu,0}^s
\rightarrow T_e \mathcal{D}_{\mu,0}^s$ be
a bounded linear map that generates a strongly-continuous
semi-group $F_t:T_e\mathcal{D}_{\mu,0}^s\rightarrow T_e\mathcal{D}_{\mu,0}^s$,
$t\ge 0$, and satisfies $\|F_t\|_s \le e^{\beta t}$ for some $\beta > 0$
and some $s$.
Extend $F_t$ to $T\mathcal{D}_{\mu,0}^s$ by
$$ \tilde{F}_t(X_\eta) = TR_\eta \cdot F_t \cdot TR_{\eta^{-1}} (X_\eta)$$
for $X_\eta \in T_\eta \mathcal{D}_{\mu,0}^s$, and let
$\tilde{{\mathcal T}}$ be the vector
field $\tilde{{\mathcal T}}:T\mathcal{D}_{\mu,0}^s \rightarrow
T^2\mathcal{D}_{\mu,0}^s$ associated to the flow $\tilde{F}_t$.

Then ${\mathcal B}+\nu \tilde{{\mathcal T}}$ generates a unique local
uniformly Lipschitz flow on
$T\mathcal{D}_{\mu,0}^s$ for $\nu \ge 0$, and the integral curves
$c^\nu(t)$ with
$c^\nu(0)=x$ extend for a fixed time $\tau >0$ independent of $\nu$ and are
unique.
Further,
$$ \lim_{\nu \rightarrow 0} c^\nu(t) = c^0(t)$$
for each $t$, $0\le t < \tau$, the limit being in the $H^s$ topology,
$s>(n/2)+1$.
\end{thm}
\begin{proof}
The proof of this theorem is essentially identical to the proof of Theorem 13.1
in \cite{EM} so will not be repeated.
\end{proof}

Now, for the equation (\ref{nse}), the operator ${\mathcal T}$ is simply the
order zero differential operator
${\mathcal T}= {\mathcal P}_e(1-\alpha^2 \triangle)^{-1} \triangle$,
coming from the equation
$$ u_t = (1-\alpha^2 \triangle)^{-1} \triangle u.$$
It is a fact that ${\mathcal T}:T_e \mathcal{D}_{\mu,0}^s \rightarrow
T_e {\mathcal D}_{\mu,0}^s$ is continuous and generates a smooth semi-group
in $T_e \mathcal{D}_{\mu,0}^s$.  This follows from the elliptic regularity
of the Stokes operator with Dirichlet boundary conditions.

Since Theorem \ref{thm_spray1}
prove that the geodesic spray ${\mathcal B}$ is $C^\infty$ on
$\mathcal{D}_{\mu,0}^s$, we use the product formula approach to iterate
the composition of the time $t/n$ maps of the vector fields ${\mathcal T}$
and ${\mathcal B}$ to obtain our result.

\begin{rem}
With initial data $u_0$ in
$T_e\mathcal{D}_{\mu,0}^\infty$, the solution $u(t)$ is
also $C^\infty$ as a consequence of the regularization of
parabolic flows.
\end{rem}

The use of the product formula in the proof of the above theorem
is given in \cite{EM}, \cite{M}, and \cite{CHMM}.

\appendix
\section{Euler-Poincar\'{e} Reduction}
\label{app1}

The reduction onto the Eulerian
representation is an example of the Euler-Poincar\'{e} theorem (see, for
example, \cite{AK1} or \cite{MR}) which we shall now state in the setting of
a topological group $G$ which is a smooth manifold and admits smooth right
translation.  For any element $\eta$ of the group, we shall denote
by $TR_\eta$ the right translation map on $TG$, so that for example, when
$G$ is either ${\mathcal N}_\mu^s$ or
$\mathcal{D}^s_{\mu,0}$, then $TR_{\eta^{-1}}{\dot{\eta}} := {\dot{\eta}}
\circ \eta^{-1}$.

\begin{thm}[Euler-Poincar\'{e}]
Let $G$ be a topological group which admits smooth manifold structure with
smooth right translation, and let $L:TG \rightarrow
{\mathbb R}$ be a right invariant Lagrangian.  Let ${\mathfrak g}$ denote
the fiber $T_eG$, and let $l:{\mathfrak g}\rightarrow{\mathbb R}$ be
the restriction of $L$ to ${\mathfrak g}$.  For a curve $\eta(t)$ in
$G$, let $u(t)=TR_{\eta(t)^{-1}}{\dot{\eta}}(t)$. Then
the following are equivalent:
\begin{itemize}
\item[\bf{a}] the curve $\eta(t)$ satisfies the Euler-Lagrange equations on
$G$;
\item[\bf{b}] the curve $\eta(t)$ is an extremum of the action function
$$S(\eta) = \int L(\eta(t),{\dot{\eta}}(t)) dt,$$
for variations $\delta \eta$ with fixed endpoints;
\item[\bf{c}] the curve $u(t)$ solves the Euler-Poincar\'{e} equations
$$ \frac{d}{dt}\frac{\delta l}{\delta u} = -\text{ad}^*_{u}
\frac{\delta l}{\delta u},$$
where the coadjoint action ad$^*_u$ is defined by
$$\langle \text{ad}^*_u v, w \rangle = \langle v, [u,w]_R\rangle,$$
for $u,v,w$ in ${\mathfrak g}$, and where $\langle \cdot, \cdot\rangle$ is
the metric on ${\mathfrak g}$ and $[ \cdot, \cdot ]_R$ is the right bracket;
\item[\bf{d}] the curve $u(t)$ is an extremum of the reduced action function
$$s(u) = \int l(u(t)) dt,$$
for variations of the form
\begin{equation}\label{con_var}
\delta u = \dot w + [w, u],
\end{equation}
where $w = TR_{\eta^{-1}}\delta \eta$ vanishes at the endpoints.
\end{itemize}
\end{thm}
See Chapter 13 in \cite{MR} for a detailed development of
the theory of Lagrangian reduction as  well as a proof of the
Euler-Poincar\'{e}
theorem.

\section*{Acknowledgments}
The authors would like to thank the anonymous referee for many helpful
suggestions that improved the manuscript.  JEM and SS were partially
supported by an NSF-KDI grant ATM-98-73133 and the DOE, and TSR was
partially supported by NSF grant DMS-98-02378 and the Swiss NSF.  SS
would also like to thank the Center for Nonlinear Studies in Los
Alamos for providing an excellent environment wherein much of this
work was performed.

\end{document}